\documentclass[11pt]{article}

\usepackage{amsmath}
\usepackage{amsfonts}
\usepackage{amsthm}
\usepackage{amsbsy}

\usepackage[english]{babel}

\usepackage{amssymb}
\usepackage{graphicx}
\usepackage{cite}
\usepackage{epstopdf}
\usepackage{epsfig}
\usepackage{subfigure}

\usepackage[usenames,dvipsnames]{color}

\usepackage{tikz}
\usetikzlibrary{shapes,arrows}
\usetikzlibrary{shapes,arrows,decorations.markings}
\usetikzlibrary{decorations.pathmorphing}

\tikzstyle{block_long} = [rectangle, draw, fill=blue!20,
    text width=10.0em, text centered, rounded corners, minimum height=3em]
\tikzstyle{block_medium} = [rectangle, draw, fill=blue!20,
    text width=6.0em, text centered, rounded corners, minimum height=3em]
\tikzstyle{block_short} = [rectangle, draw, fill=blue!20,
    text width=3.0em, text centered, rounded corners, minimum height=3em]
\tikzstyle{line} = [thick, draw, dashed,  -stealth']

\tikzstyle{block_long} = [rectangle, draw, fill=blue!20,
    text width=13.0em, text centered, rounded corners, minimum height=3.5em]
\tikzstyle{block_medium} = [rectangle, draw, fill=blue!20,
    text width=7.0em, text centered, rounded corners, minimum height=3.5em]
\tikzstyle{block_short} = [rectangle, draw, fill=blue!20,
    text width=5.0em, text centered, rounded corners, minimum height=3.5em]
\tikzstyle{line} = [thick, draw, dashed,  -stealth']
\tikzstyle{chapter} = [rectangle, draw, fill=blue!10,
text width=6.0em, text centered, rounded corners, minimum height=0.5em]

\usepackage{multirow}
\usepackage{amsmath}
\usepackage{amssymb}
\usepackage{soul}
\usepackage{ulem}

\usepackage[mathscr]{eucal}

\newtheorem{Lemma}{Lemma}

\usepackage{color}

\def\bw{{\bf w}}

\def\bu{{\bf u}}
\def\bx{{\bf x}}

\def\vf{{\varphi}}
\def\ve{{\varepsilon}}

\def\bof{{\bf f}}

\def\div{\mbox{{\rm div}}}

\def\va{\raise 2pt\hbox{,}}

\def\@xthm#1#2{\@beginassumption{#2}{\csname the#1\endcsname}{}\ignorespaces}
\def\@ythm#1#2[#3]{\@opargbeginassumption{#2}{\csname the#1\endcsname}{#3}\ignorespaces}%
\def\@beginassumption#1#2#3{\par\addvspace{8pt plus3pt minus2pt}%
              \noindent{\csname#1headfont\endcsname#1\ \ignorespaces#3 #2.}%
              \csname#1font\endcsname\hskip.5em\ignorespaces}
\def\@endassumption{\par\addvspace{8pt plus3pt minus2pt}\@endparenv}
\title{Cross diffusion models in complex frameworks\\[3mm] From microscopic to macroscopic.}

\author{
D.~Burini
\thanks{Universit\`a degli Studi di Perugia, Perugia, Italy.
 {(dilettaburini@alice.it)}
}
\and
N.~Chouhad
	\thanks{ Cadi Ayyad University, Ecole Nationale des Sciences Appliqu\'ees, Marrakech, Morocco,
{(chouhadn@gmail.com)}
			}
}

\begin{document}
\maketitle

\begin{abstract}
This paper deals with the micro-macro derivation of models from the underlying description provided by methods of the kinetic theory for active particles. We consider the so-called \textit{exotic} models according to the definition proposed in in\cite{[BOSTW22]}. The first part  of the  presentation focuses on a survey and a critical analysis of some phenomenological models known in the literature. We refer to a selection of case studies, in detail: transport of virus models, social dynamics, and Keller-Segel in a fluid. The second part shows how an Hilbert type approach can be developed to derive  models at the macroscale from the underlying description provided by the kinetic theory of active particles. The third part deals with the derivation of macroscopic models corresponding to the selected case studies. Finally, a forward look into the future research perspectives is proposed.
\vskip.2cm


\vskip.2cm


\end{abstract}

\section{Aims and plan of the paper}\label{Sec:1}

A review and critical analysis, on the modeling and analytical problems of the so-called \textit{exotic cross diffusion systems}, was recently developed in\cite{[BOSTW22]}. These models describe the dynamics of living systems within the general framework of cross-diffusion and reaction-diffusion models,  for short CD and RD models\cite{[KS70],[KS71]}. Different types of models can be, so-far, defined \textit{exotic}. For instance,  models whose dependent variable corresponds to social systems, e.g., modeling social dynamics\cite{[SBBT10]},  models driving a virus\cite{[BPTW19],[BT20]}, and CD and RD models within a fluid\cite{[BBC16]}. Additional models can be obtained by possibly mixing of the various features of the above examples.

The derivation of models  in\cite{[BOSTW22]} is essentially related to the phenomenological approach corresponding to conservation/transport equations closed by heuristic models considered to describe the material behavior of the system. The survey of analytical studies generally correspond to initial-boundary value problems. Some hints towards a further refinement and improvement are therein reported, but the derivation of models at the macroscopic scale from the underlying description at the microscopic scale is still open with the exception of the preliminary results on the transport of virus models by a Keller-Segel system\cite{[BC22]}, briefly \textit{KS system}. The contents of  previous papers  by the same authors\cite{[BC17],[BC19]} will contribute to the development of the present article.

Our  paper refers to the micro-macro derivation of exotic models from the underlying description provided by methods of the kinetic theory of active particles\cite{[BBDG21]} by further developments of the  method developed in\cite{[BC22]} somehow inspired by the celebrated sixth Hilbert problem\cite{[Hilbert],[SLEM13]}, generally related to the kinetic theory of classical particles\cite{[GK14],[SLEM13]}, but subsequently developed for biological systems\cite{[BC17]}. We focus mainly  on the micro-macro derivation of virus models with space dynamics. However, possible developments concerning different types of applications are also brought to the attention of the interested reader. The contents are as follows:

\vskip.1cm Section 2 presents an overview of exotic models with reference to a selection of models. Specifically, we consider different types  of virus models and show that the space dynamics can be heuristically modeled by cross reaction-diffusion models with a source term. We do start from\cite{[BOSTW22]} and then we present various developments concerning both the structure of the CD and RD models and the source term.  This section also considers some classes of social models. A general framework at the macroscopic scale is derived to encompass all the models presented in this section.

 \vskip.1cm Section 3 deals with the first step  towards the derivation of macroscale models from the underlying microscale description. In detail,  a general methodological approach  is derived based on developments, with respect to those proposed in\cite{[BC22]},  of the Hilbert approach. We refer to the general structure proposed in Section 2.

\vskip.1cm Section 4  shows how the approach  can be applied to the micro-macro derivation of the specific  models presented  in Section 2. These applications concern both  virus models and to social systems.

\vskip.1cm Section 5 searches for new classes of models, still within a framework analogous to that of Section 2, but with source terms  modeled by the kinetic theory of active particles.

\section{Exotic cross diffusion models}\label{Sec:2}

This section reports on exotic cross-diffusion models through a selection of case studies to be considered for the development of a micro-macro derivation by a Hilbert-type approach. The focus  is on the spatial dynamics of biological and social models, described by systems of ODEs,  acting as a source term of a  reaction-diffusion system.

First, a general framework is presented with the idea of incorporating a wide variety of case studies. A selection of specific models is then presented. In detail, we first briefly review  the population dynamics modeling of virus dynamics under the action of a reaction-diffusion system and we propose the aforementioned  mathematical structure.  We then report some case studies that can be derived within such a specific framework. These case studies selected as applications for the micro-macro derivation developed in Sections 3 and 4.

\subsection{ODEs models under the action of a reaction-diffusion system}\label{Sec:2.1}

Let us consider mathematical models derived, within the general framework of population dynamics, derived within the general framework of vector n-dimensional ODEs:
\begin{equation}\label{ODEs}
\frac{d \bu}{dt} = \bof (t, \bu), \hskip1cm \bu =\{u_1, \ldots, u_n\}, \hskip1cm \bof =\{f_1, \ldots, f_n\},
\end{equation}
where $\bu$ is an n-dimensional variable modeling the state of a system of $n$ interacting populations labeled by the subscript $i= 1, \ldots, n$.
A classical example is the May-Nowak model which describes the dynamics of three components, i.e. the densities of healthy uninfected immune cells $u_1= u_1(t)$, infected immune cells $u_2= u_2(t)$, and virus particles $u_3 = u_3(t)$, see~\cite{[NB96],[NM00]}  for the assumptions leading to derivation of the model, see also for the study of the qualitative properties of the model~\cite{[PERT07]}.

Different variation of models are known in the literature, starting from~\cite{[NM00]}, where the study of specific aspects of the immune competition are treated referring to biological applications. Several developments  followed also in recent years mainly to investigate specific aspects of the immune competition. Generally, models include a greater number of populations and delay terms, as examples, see~\cite{[EA17],[LM14]}. However, the aims of our paper do not include a detailed report on the literature on population dynamics, we simply address the interested reader to some books devoted to the mathematical and biological study of this topic~\cite{[DH00],[NOW06],[NM00],[PERT07]}.

These, rapidly reviewed, models do not describe dynamics in space. However, various authors have studied models, where $\bu = \bu(t, \bx)$ depends on space $\bx$ in addition to time $t$, and where space dynamics is modeled by a reaction-diffusion system acting on the variable $\bu$. As a particular case, the reaction-diffusion action term corresponds to a  Keller-Segel chemotaxis system~\cite{[BPTW19]}. Specific applications have been studied in in~\cite{[BPTW19],[BT20],[BMSN97],[FUEST19]}, while a qualitative of mathematical problems has been developed by various authors, as reviewed in~\cite{[BOSTW22]}.

The mathematical structure in the following is modeled by a deterministic reaction diffusion dynamics acting on  $u_1 = u_1(t,\bx)$, $u_2 = u_2(t,\bx)$ and $u_3 = u_3(t,\bx)$ which now include  space dependence, just as a tutorial example, to antucuoate how more sophisticated models can be referred to this simple case study:
\begin{equation}\label{1}
\begin{cases}
  \displaystyle \partial_t u_1 = D_1 \Delta u_1 - \chi \nabla\cdot (u_1\nabla u_2) - d_1 \, u_1 - \beta \, u_1 \, u_3 + r(t, \bx), \\[4mm]
  \displaystyle  \partial_t u_2 = D_2 \Delta u_2 - d_2 \,u_2 + \beta \,u_1 \, u_3,   \\[4mm]
    \displaystyle \partial_t u_3 = D_3 \Delta u_3 - d_3\, u_3 + k\, u_2,  \\
     \end{cases}
    \end{equation}
where $D_1$, $D_2$ and $D_3$ denote the respective, positive defined, diffusion coefficients and where $\chi$ represents strength and direction of the cross-diffusive interaction. Healthy cells are constantly produced by the body at rate $r$, die at rate $d_1 u_1$ and become infected on contact with the virus, at rate $\beta u_1  u_3$; infected cells are  produced at rate $\beta u_1  u_3$ and die at rate $d_2 u_2$, while new virus particles are produced at rate $k u_2$ and die at rate $d_3 u_3$.  The reaction-diffusion action term corresponds to a simplified  Keller-Segel chemotaxis system~\cite{[BPTW19],[BT20]}.

The models presented in the following  refer to the following general structure:
\begin{equation}\label{virus-a}
\begin{cases}
\displaystyle{\partial_t u_1 = \vf_1\left(\bu, \Delta \bu, \nabla \bu; \bw \right) + \psi_1(\bu)},\\
\displaystyle{ \ldots ,}\\[2mm]
\displaystyle{\partial_t u_i = \vf_i\left(\bu, \Delta \bu, \nabla \bu; \bw \right) + \psi_i(\bu)},\\
\displaystyle{ \ldots ,}\\[2mm]
\displaystyle{\partial_t u_n = \vf_n\left(\bu, \Delta \bu, \nabla \bu; \bw \right) + \psi_n(\bu)}. \\
\end{cases}
\end{equation}

Various examples are presented in the next subsections in view of the micro-macro derivation. In some cases, the dynamics of one or more components of $\bu$ do not include space diffusion, so that the mathematical structure of these components is simply an ODE.

\subsection{Oncolytic viruses}\label{Sec:2.2}

An example which refers to the framework presented in Subsection 2.1   is the model  that that describes  selective replication of a virus within cancer cells and attack them up to destruction, see~\cite{[AET20],[TLCW13]}. The model consists in a 4-dimensional system coupling partial differential equations at the macro-scale (tissue-scale) with a  system at the micro-scale (cell-scale). In more detail, the model describes the dynamic interactions between four macro-scale components moving on a domain $\Omega(t)$, which evolves with time. These components are the density of uninfected cancer cells $u_1= u_1(t, \bx)$, the density of infected cancer cells $u_2= u_2(t, \bx)$, the density of extracellular matrix $u_3= u_3(t, \bx)$, and the density of the oncolytic virus particles $u_4= u_4(t, \bx)$, with $\bx \in \Omega(t)$ and $t>0$.

In its general form, the model reads as follows:
\begin{equation}\label{oncolytic}
\begin{cases}
\displaystyle{ \partial_t u_1 = D_1 \Delta u_1 - \xi_1 \nabla \cdot \left(u_1\nabla u_3\right) + \mu_1 u_1 (1 - u_1^r) - \frac{\rho u_1 u_4}{k_{u_1}+\theta u_1}},\\[4mm]
\displaystyle{ \partial_t u_2 = D_2 \Delta u_2 - \xi_2 \nabla \cdot \left(u_2 \nabla  u_3\right) + \frac{\rho u_1 u_4}{k_{u_1}+\theta u_1} - \delta_2 u_2,}\\[4mm]
\displaystyle{ \partial_t u_3 = - u_3\left( \alpha_1 u_1 +  \alpha_2 u_2\right) + \mu_3 u_3(1- u_3),}\\[4mm]
\displaystyle{ \partial_t u_4 = D_4 \Delta u_4  - \xi_4 \nabla \cdot \left(u_4 \nabla u_3\right) + \beta u_2-\delta_4 u_4 - \frac{\rho u_1 u_4}{k_{u_1}+\theta u_1}}.
\end{cases}
\end{equation}

Let us observe that, in addition to random diffusion with the respective motility coefficient $D_1$ and $D_2$, the model describes the dynamics of cancer cells which direct their movement toward regions of higher extra cellular matrix, shortly ECM, densities with the haptotactic coefficients $\xi_1$, $\xi_2$, respectively, and that uninfected cells, apart from proliferating logistically at rate $\mu_1$, are converted into an infected state upon contact with virus particles.

The ECM is supposed  to be continuously remodeled by cells in the environment. This remodeling process is modeled as the difference between a logistic growth term (describing the deposition of ECM components in the presence of cancer  at a rate $\mu_3$) and a degradation term (with $\alpha_1$ the rate of ECM degradation by uninfected cancer cells, and $\alpha_2$ the rate of ECM degradation by infected cancer cells).

It is also supposed that besides the random motion with $D_4$ the random motility coefficient, virus particles move up the gradient of ECM
with the ECM-OV-taxis rate $\xi_4$, increase at a rate $\beta$ due to the release of free virus particles through infected cells and undergo decay at the rate $\delta_4$ accounting for the natural virions' death as well as the trapping of these virus particles into the cancer
cells.

In (\ref{oncolytic}) the haptotactic motion of virus particles is taken into account particularly, and the production term $u_1 u_4$
has been replaced by
$$
\frac{\rho u_1 u_4}{k_{u_1}+\theta u_1}
$$
corresponding to the Beddington--deAngelis type model with positive parameters $k_{u_1}$, $\theta$ while the proliferating term $\mu_1 u_1 (1-u_1)$ is adjusted to $\mu_1 u_1 (1-u_1^r)$  with a positive parameter $r$. For more details see \cite{[KLW22]}, Chapter 6, Section 6.1, page 347.

\subsection{In-host biological dynamics related to cancer modeling}\label{Sec:2.3}

Reaction-diffusion equations have been developed to model pattern formation in cancer biology. The physical-biological reality is highly complex as it includes mutations and selection followed by proliferative and/or destructive events related to  the interaction between tumor cells and immune cells~\cite{[HW11]}, see also~\cite{[MF19]} for the biology of the immune competition, as well as  some pioneering interpretations developed in~\cite{[BD06]}, see also~\cite{[AET20],[TLCW13]}.
Our subsection, which is focused on modeling, accounts for the analytic results in~\cite{[PW18],[PW20]}, of the surveys~\cite{[CL06],[WANG20]}, see also~\cite{[KLW22]}. The model writes
\begin{equation}\label{Lolas}
\begin{cases}
  \displaystyle {\partial_t  u_1 =  \Delta u_1 - \chi \nabla\cdot (u_1\nabla u_2) - \xi  \nabla\cdot (u_1\nabla u_3) + \mu\, u(r -u_1 - u_3)}, \\[4mm]
  \displaystyle {\partial_t u_2  = \frac{1}{\sigma} \Delta u_1 -  \frac{1}{\sigma} (u_2 - u_1)\va }  \\[4mm]
    \displaystyle {\partial_t  u_3 = - u_2 \,  u_3 + \eta\,  u_3(1 - u_1 -  u_3)},
     \end{cases}
    \end{equation}
where $u_1 = u_1(t,\bx)$,  $u_2 = u_2(t,\bx)$, and  $ u_3 =  u_3(t,\bx)$ are the model variables corresponding to the density of cancer cells, the concentration of the Matrix-Degrading Enzyme (MDE), and  the concentration of the Extracellular Matrix (ECM), respectively.
 The independent variables are time $t$ and space  $\bx \in \Omega$, where $\Omega$ is a bounded domain with regular surface  $\partial\, \Omega$  so that the outward normal derivative on $\partial\, \Omega$ can be defined for the statement of the boundary value problems.

The biological meaning of the parameters and source terms are as follows:

\vskip.1cm \noindent  $\chi$  models the chemotactic sensitivity;

\vskip.1cm \noindent  $\xi$ models the haptotactic sensitivity;

\vskip.1cm \noindent   $\mu(r - u_1 -  u_3)$ implies that in the absence of the ECM, cancer cells proliferate according to
 a standard logistic law;

\vskip.1cm \noindent  $\eta > 0$ models the ability of the ECM to remodel back to a healthy level as a coefficient to the term $u_3(1 - u_1 -  u_3)$.

\vskip.1cm   The scaling parameter $\sigma \in [0,1]$, where the limit values $\sigma = 0$ and  $\sigma = 0$, define, respectively, two biological limit behaviors; specifically,  $\sigma = 0$ indicates that the diffusion of the enzyme is much faster compared to that of cancer cells~\cite{[CL06]}, which may also follow an approach of the quasi-steady-state approximation frequently used to study minimal chemotaxis systems.

\subsection{Cross-diffusion of social (biological) dynamics}\label{Sec:2.4}

As an example of social dynamics, we consider models, where specific populations compete to chase the same objective, for instance, food distributed in space.  These populations undergo a cross-diffusion dynamics somehow promoted by such objective acting as a source, see~\cite{[TVHE12]}.
 Analytic topics concerning the qualitative analysis of solutions to  initial-boundary value problems are reviewed in~\cite{[BOSTW22]}.

In details, we consider models where the space dynamics  is modeled by reaction-diffusion equations chasing a source that depends on the space availability of food, which diffuses being somehow reduced by the action of the two populations. Specifically, \textit{foragers} that search for food directly,  and \textit{scroungers}, (say \textit{exploiters}), in the search the food by following the foragers, namely exploiting  forager aggregations. The  dynamics of food interacts with that of the first two populations. A quite general model is as follows:
\begin{equation}\label{fora-exploit}
\begin{cases}
\displaystyle{\partial_t u_1 =  \Delta \,u_1 - \xi_1  \nabla \cdot \left(u_1 \nabla \, u_3 \right) + f(u_1, u_2, u_3)},\\[3mm]
\displaystyle{\partial_t u_2 =  \Delta \,u_2 - \xi_2  \nabla \cdot \left(u_2 \nabla \, u_1 \right) + g(u_1, u_2, u_3)},\\[3mm]
\displaystyle{\partial_t u_3 = D \Delta \,u_3 - \lambda(u_1 + u_2)\,u_3 - \mu\,u_3 +r(u_3)},
\end{cases}
\end{equation}
where $u_1=u_1(t,\bx)$, $u_2=u_2(t,\bx)$, and $u_3=u_3(t,\bx)$ are, respectively, the dimensionless densities of foragers, scroungers, and food, while $f(u_1, u_2, u_3)$ and $g(u_1, u_2, u_3)$ are  source terms acting across the first two populations, and $r(u_3)$ is the inner source of the third population. Dimensionless variables are used, while the  dimensionless parameters $\xi_1$ and $\xi_2$ are  taxis parameters of foragers and exploiters, $D$ is the relative mobility of the prey, $\lambda$ is the per-capita consumption rate,  and $\mu$ is the decay rate of the prey. This framework defines the so-called \textit{taxis-cascade systems}, see~\cite{[KLW22]}.

\section{An Hilbert method towards derivation of macro-scale models}\label{Sec:3}

This section presents a theoretical approach towards the derivation  macro-scale models from the underlying description at the micro-scale delivered within a kinetic theory framework. We consider models which can be derived according to the general structure given by Eq.~(2.3). The contents is in two steps developed in the next two subsections. Firstly we propose a kinetic theory structure obtained as perturbation of a transport term. Then we show how a general macroscopic model can be obtained by an Hilbert-type expansion.
Application to the derivation of specific models will be developed in the next section focusing on the case studies presented in Section 2.

\subsection{On a kinetic theory framework towards the asymptotic analysis}\label{Sec:3.1}

This subsection presents a general mathematical structure corresponding to a micro-macro decomposition for a n-components mixture of self-propelled particles  whose state, called {\sl microscopic state}, is denoted by the variable $(x, v)$, corresponding, respectively, to position and velocity.  The collective description  of a mixture of particles  can be encoded in the statistical distribution functions $f_i = f_i(t, x, v)$, for $i= 1, 2, \cdots, n$. Weighted moments provide, under suitable integrability properties, the calculation of macroscopic variables.

Let us now consider the following class of equations:
\begin{equation}\label{LinearT}
\begin{cases}
 \big(\partial_t + v \cdot \nabla_{x} \big) f_1 =\nu_1\, \mathcal{T}_1[f_1,f_2,\cdots,f_n](f_1)+\mu_1\,G_1(f_1,f_2,\cdots,f_n,v),\\[4mm]
   \big(\partial_t + v \cdot \nabla_{x}\big)f_2=\nu_2\, \mathcal{T}_2[f_1,f_2,\cdots,f_n](f_2 )+ \mu_2\, G_2(f_1,f_2,\cdots,f_n,v), \\[4mm]
   \big(\partial_t + v \cdot \nabla_{x} \big) f_3= \nu_3\, \mathcal{T}_3(f_3 ) +\mu_3\,G_3(f_1,f_2,\cdots,f_n,v),     \\[4mm]
  \big(\partial_t + v \cdot \nabla_{x}\big)f_4=\nu_4\, \mathcal{T}_4[f_1,f_2,\cdots,f_n](f_4 )+ \mu_4\, G_4(f_1,f_2,\cdots,f_n,v),\\[4mm]
\cdots, \\[4mm]
 \big(\partial_t + v \cdot \nabla_{x} \big) f_n= \nu_n\,\mathcal{T}_n[f_1,f_2,\cdots,f_n](f_n )+ \mu_n \,G_n(f_1,f_2,\cdots,,f_n,v),
\end{cases}
\end{equation}
where  $G_1, G_2 , \cdots, G_n$ are interactions terms assumed depending on the quantities
$f_1, f_2, \cdots, f_n$, while the operator $\mathcal{T}_i(f)$ models the dynamics of biological organisms by a velocity-jump process:
\begin{equation}\label{de}
\mathcal{T}_i(f)= \int_{V} \bigg[T_i(v^*, v)f(t, x, v^{*}) - T_i(v , v^*)f(t, x, v) \bigg]\,dv^{*}, \quad i=1,2,\cdots ,n,
\end{equation}
where $T_i(v, v^{*})$ is the probability kernel for the new velocity $v \in V$ assuming that the previous velocity was $v^{*}$.

System (\ref{LinearT}) can be put in a dimensionless form, see \cite{[BC19]}, so that a small parameter $\varepsilon$ can be extracted corresponding to the parabolic scaling:
\[t\longrightarrow \varepsilon t, \quad \mu_i=\varepsilon, \quad
\nu_\ell=\frac{1}{\varepsilon},\quad  \nu_3=\frac{1}{\varepsilon^q},   ~~~i=1,2,\cdots,n,~~~ \ell=1,2,4,\cdots,n,~~~q>1.
\]
Then, the model (\ref{LinearT}) can be rewritten as  follows:
\begin{equation}\label{kinetic}
\left\{
\begin{array}{l}
\big(\varepsilon \partial_t + v \cdot \nabla_{x} \big) f_\ell^\varepsilon= \displaystyle{\frac{1}{\varepsilon}} \mathcal{T}_\ell[f^\varepsilon_1,f^\varepsilon_2,\cdots,f^\varepsilon_n](f_\ell^\varepsilon ) +\varepsilon\,G_\ell(f^\varepsilon_1,f^\varepsilon_2,\cdots,f^\varepsilon_n, v), \\[4mm]
\big(\varepsilon \partial_t + v \cdot \nabla_{x} \big) f_3^\varepsilon  =  \displaystyle {\frac{1}{\varepsilon^q} \mathcal{T}_3(f_3^\varepsilon )} +\varepsilon\,G_3(f^\varepsilon_1,f^\varepsilon_2,\cdots, f^\varepsilon_n, v).
\end{array}\right.
\end{equation}

The  derivation of macroscopic models from the kinetic model (\ref{LinearT}), can be obtained by letting $\varepsilon \, \to \,0$. The singular perturbation requires appropriate assumptions  on the terms $\mathcal{T}$ modeling the stochastic perturbation. These assumptions are  reported in the following.

\vskip.2cm \noindent {\bf Assumption 3.1.} \label{assH01} The   turning operators  $\mathcal{T}_1, \mathcal{T}_2, \cdots $ , $\mathcal{T}_n$, are supposed to be decomposable as follows:
 \begin{equation}
\mathcal{T}_{\ell}[ f^\varepsilon_1,f^\varepsilon_2,\cdots,f^\varepsilon_n](g)= \mathcal{T}_{\ell}^{0}(g)+\varepsilon^{b_{\ell}} \, \mathcal{T}_{\ell}^{1}[ f^\varepsilon_1,f^\varepsilon_2,\cdots,f^\varepsilon_n](g), \quad  b_{\ell} \geq 1,\quad  \ell =1,2,4,\cdots,n,\label{L1}
\end{equation}
where $\mathcal{T}_{\ell}^{j}$ for $j=0,1$, is given by

\begin{equation}
\mathcal{T}_{\ell}^{j}(g)= \int_{V} \bigg[T_{\ell}^{j}(v^{*},v)g(t, x, v^*) - T_{\ell}^{j}(v,v^{*}) g(t, x, v) \bigg]dv^{*},  \label{3.4}
\end{equation}

 where the dependence on  $f_1,\cdots, f_n$, of the operator $\mathcal{T}_{\ell}$ stems from $T_{\ell}^1$, while we suppose that $\mathcal{T}_{\ell}^0$ is independent of $f_1,\cdots, f_n$,  and  $\mathcal{T}_{3}$ is
\begin{equation}
\mathcal{T}_{3}(g)= \int_{V} \bigg[T_{3}(v^{*},v)\,g(t, x, v^{*}) - T_{3}(v,v^{*}) g(t, x, v) \bigg]dv^{*}.
\end{equation}

\vskip.2cm \noindent {\bf Assumption 3.2.}  We assume that the turning operators $\mathcal{T}_i (i=1,2,\cdots,n)$  satisfy the following equality:
\begin{equation}\label{assH3}\int_V \mathcal{T}_i(g)dv=\int_V\mathcal{T}_\ell^0(g)dv= \int_V \mathcal{T}_\ell^1[f^\varepsilon_1,f^\varepsilon_2,\cdots,f^\varepsilon_n](g)dv=0,\quad  \ell =1,2,4,\cdots, n.
\end{equation}

\vskip.2cm \noindent {\bf Assumption 3.3}\label{assH4} There exists a bounded velocity distribution $ M_3(v)>0$ and  $M_\ell(v)>0$ for $\ell =1,2,4,\cdots,n$, independent of $t, x$, such that the detailed balance
\begin{equation}
T_\ell^0 (v,v^{*} ) M_\ell(v^{*}) = T_\ell^0 (v^{*},v ) M_\ell(v),
\end{equation}
and
\begin{equation}\label{tx1}
T_3(v,v^{*} ) M_3(v^{*}) = T_3(v^{*},v ) M_3(v),
\end{equation}
hold true. Moreover, the flow produced by these equilibrium distributions vanishes, and $M_i$ are normalized
\begin{equation}\label{TZ}
\int_V v \, M_i(v)dv  =0, \quad \int_V
M_i(v)dv = 1,\quad i=1,2,\cdots,n.
\end{equation}

\vskip.2cm The kernels  $T_3(v,v^{*})$ and  $T_\ell^0(v,v^{*})$  are bounded and that there exist constants  $\sigma_3>0$ and $\sigma_\ell>0$, $\ell=1,2,4,\cdots,n$, such that
\begin{eqnarray}\label{cx}
T_3(v,v^{*})\geq \sigma_3 M_3(v),\quad T_\ell^0(v,v^{*})\geq \sigma_\ell M_\ell(v),
\end{eqnarray}
for all $ (v,v^{*}) \in V\times V$, $x\in \Omega$ and $t > 0$.
\vskip.2cm

Given that $L_3=\mathcal{T}_3$ and $L_\ell=\mathcal{T}^0_\ell(\ell=1,2,4,\cdots,n)$. Technical calculations yields the following Lemma:
\vskip.1cm
\begin{Lemma}\label{LE1}
Suppose that  Assumption 3.3  holds. Then, for $i=1,2,\cdots,n$ the following properties of the operators $L_i$ hold:
\begin{itemize}
\item[i)] The operator $L_i$ is self-adjoint in the space $\displaystyle{{L^{2}\left(V ,{dv \over M_i}\right)}}$.
\vskip.1cm \item[ii)]
For $f\in L^2$, the equation $L_i(g)=f$ has a unique solution $\displaystyle{g \in L^{2}\left(V, {dv \over M_i}\right)}$, which satisfies
$$
 \int_{V} g(v)\, dv = 0, \quad  \hbox{if and only if} \quad   \int_{V} f(v)\, dv =0.
$$
\vskip.1cm \item[iii)]  The equation $L_i(g) =v \,  M_i(v)$ has a unique solution that we call $\theta_i(v)$.
\vskip.1cm \item[iv)]  The kernel of $L_i$ is $N(L_i) = vect(M_i(v))$.
\end{itemize}
\end{Lemma}

\subsection{Derivation of a general  macroscopic models} \label{Sec:3.2}

A system coupling a hydrodynamic part with a kinetic part of the distribution functions, is derived in this subsection. Then it is proved that  such a system is equivalent to the two scale kinetic equation (\ref{kinetic}). This new formulation provides the basis for the derivation of the  general model we are looking for.

In the remainder, the integral with respect to the variable $v$ will be denoted by $\langle \cdot  \rangle$. This notation is used also for  any argument  within $\langle \, \rangle$. In addition,   let us denote by  $f=(f_1, f_2,\cdots,f_n)$ the solution of (\ref{kinetic}), where $f$  is decomposed as follows:
\begin{equation} \label{eq}
f_{\ell}^\varepsilon(t,x,v)=\sum^{2}_{k=0}\varepsilon^{k}g_{ k\ell}(t,x,v)+O(\varepsilon^{3}),\quad \ell=1,2,4,\cdots,n,
  \end{equation}
and
\begin{equation} \label{eq1}
f_{3}^\varepsilon(t,x,v)=\sum^{q+1}_{j=0}\varepsilon^{j}h_{j}(t,x,v)+O(\varepsilon^{q+2}).
 \end{equation}

In order to develop asymptotic analysis of Eq.~(\ref{kinetic}),  additional assumptions on the operator $\mathcal{T}_\ell^1(\ell=1,2,4,\cdots,n)$  and the interaction terms $G_{i} (i=1,2,\cdots,n)$   are needed.

\vskip.2cm \noindent {\bf Assumption 3.4.}
We assume that the turning operator $\mathcal{T}_\ell^1(\ell=1,2,4\cdots,n)$  and the interaction terms $G_{i}(i=1,2,3,\cdots,n)$
 satisfies the following asymptotic behavior as:
\begin{equation}\label{R1}
\begin{cases}
\mathcal{T}_\ell^1\big[f^\varepsilon_1,f^\varepsilon_2,\cdots,f^\varepsilon_n](g)= \mathcal{T}_\ell^1[g_{01},g_{02},h_{0},g_{04},\cdots,g_{0n}] (g)\\[4mm]
\hskip1cm + \varepsilon\,\mathcal{R}_\ell^1[ g_{01},g_{11},g_{21},\cdots, h_0,\cdots ,h_q,\cdots ,g_{0n},g_{1n},g_{2n}](g)\\[4mm]
\hskip1cm + O(\varepsilon^{2}) ,~~ \forall b_\ell=1,~~ \forall g,\\[4mm]
\mathcal{T}_\ell^1\big[f^\varepsilon_1,f^\varepsilon_2,\cdots,f^\varepsilon_n](g)= \mathcal{T}_\ell^1[g_{01},g_{02},h_{0},g_{04},\cdots,g_{0n}] (g)+ O(\varepsilon) ,~ \forall b_\ell> 1,~ \forall g.
\end{cases}
\end{equation}
and
\begin{equation}\label{G21}
  G_{i}(g_1+\varepsilon \,  \hat{g_1}, g_2+\varepsilon \,  \hat{g_2}, ....,g_n+\varepsilon \,  \hat{g_n},v) =G_{i}(g_1,g_2,....,g_n,v)+ O(\varepsilon), \hskip.2cm \forall g_1,\hat{g_1},..., g_n,\hat{g_n}.
 \end{equation}

Then, from (\ref{assH3}) one has:
\begin{equation}\label{R12}
\int_{V} \mathcal{R}_\ell^1[ g_{01},g_{11},g_{21},\cdots, h_0,\cdots ,h_q,\cdots ,g_{0n},g_{1n},g_{2n}](\varphi)\,dv=0,  \quad \forall \varphi.
 \end{equation}

\vskip.2cm Therefore, the first terms of Hilbert expansion of equal order in $\varepsilon^k$ and  $\varepsilon^j$ for $k=0,1,2$, $j=0,1,2, \cdots, q+1$ and  $ \ell =1,2,4,\cdots, n$  are:
\begin{equation}\label{c1}
\ve^0:
 \left\{
 \begin{array}{l}
 \mathcal T_{\ell}^{0}(g_{0\ell})=0, \\[4mm]
\mathcal T_{3}(h_{0})=0,
\end{array}  \right.
\end{equation}

\begin{equation}\label{ccoo}
\ve^1:
 \left\{
 \begin{array}{l}
 \mathcal T_{\ell }^{0}(g_{1\ell })= v\cdot \nabla_{x} g_{0\ell }- \delta_{b_\ell ,1}\mathcal T_{\ell }^{1}[g_{01},g_{02},h_{0},g_{04},\cdots,g_{0n}](g_{0\ell }),\\[4mm]
 \mathcal T_{3}(h_{1})=\delta_{q,1}\, v\cdot \nabla_{x} h_{0},
\end{array}  \right.
\end{equation}

and
\begin{equation}\label{ccc}
{\ve^2:}
 \left\{
 \begin{array}{llll}
 \mathcal T_{\ell }^{0}(g_{2\ell })= \partial_t g_{0\ell }+ v\cdot \nabla_{x} g_{1\ell }-\delta_{b_\ell ,2} \mathcal T_{\ell }^{1}[g_{01},g_{02},h_{0},g_{04},\cdots,g_{0n}](g_{0\ell }) \\[4mm] \hskip0.5cm-\delta_{b_\ell ,1} \mathcal T_{\ell }^{1}[g_{01},g_{02},h_{0},g_{04},\cdots,g_{0n}](g_{1\ell })  \\[4mm]
\hskip0.5cm-\delta_{b_\ell,1} \mathcal R_{\ell }^{1}[ g_{01},g_{11},g_{21},\cdots, h_0,\cdots ,h_q,\cdots ,g_{0n},g_{1n},g_{2n}](g_{0\ell })\\[4mm]
\hskip0.5cm- G_{\ell }(g_{01},g_{02} ,h_{0},g_{04},\cdots, g_{0n} ,v),\\[4mm]
\mathcal T_{3}(h_{2})= \delta_{q,1} \big(\partial_t h_{0} + v\cdot \nabla_{x} h_1 \big) + \delta_{q,2}\, v\cdot \nabla_{x}{h_0}\\[4mm]
 \hskip0.5cm -\delta_{q,1}G_{3}(g_{01},g_{02} ,h_{0},g_{04},  \cdots, g_{0n} ,v),
\end{array}
 \right.
\end{equation}
Further calculations yield:
\begin{equation} \label{ccc2}
\varepsilon^{q+1}:\mathcal T_{3}(h_{q+1})= \partial_t h_{0} + v\cdot \nabla_{x} h_1 -G_{3}(g_{01},g_{02} ,h_{0},g_{04},\cdots, g_{0n} ,v),
\end{equation}
where $\delta_{a,b}$ stands for the Kronecker delta.
The first equation of (\ref{c1}) implies that
$$
g_{0\ell}\in vect(M_{\ell}(v)),\quad   h_{0}\in  vect(M_{3}(v)).
$$

Therefore $\exists\,u_{\ell}(t,x)(\ell=1,2,4,\cdots n),$ and $\exists\,u_{3}(t,x)$ such that
\begin{equation} \label{nS1}
g_{0\ell}(t,x,v)= M_\ell(v)\,u_\ell(t,x),  ~~h_0(t,x,v)= M_3(v)\,u_3(t,x).
\end{equation}
Using (\ref{assH3}), (\ref{TZ}) and  (\ref{nS1}) we conclude that Eq. (\ref{ccoo}) satisfies the solvability condition, therefore $g_{1\ell},$ and  $h_{1}$  are given by
\begin{equation}\label{ph}
  \begin{cases}
g_{1\ell}=(\mathcal T_{\ell}^{0})^{-1}(v\cdot\nabla_{x} g_{0\ell})- \delta_{b_\ell,1}(\mathcal T_{\ell}^{0})^{-1}(\mathcal T_{\ell}^{1}[g_{01},g_{02},h_{0},g_{04},\cdots,g_{0n}](g_{0\ell})),
 \\[4mm]
 h_{1}=\delta_{q,1} \mathcal T_{3}^{-1}(v\cdot\nabla_{x} h_{0}),
\end{cases}
\end{equation}
The calculations of $g_{2\ell}$, and  $h_{q+1}$  are obtained from the solvability conditions at $O(\varepsilon^{2})$ and $ O(\varepsilon^{q+1})$, which are given by the following:
\begin{equation}\label{cc1}
 \left\{
 \begin{array}{llll}
\displaystyle\int_{V}\bigg( \partial_t g_{0\ell }+ v\cdot \nabla_{x} g_{1\ell }-\delta_{b_\ell ,2} \mathcal T_{\ell }^{1}[g_{01},g_{02},h_{0},g_{04},\cdots,g_{0n}](g_{0\ell })\\[4mm]
\hskip0.5cm-\delta_{b_\ell ,1} \mathcal T_{\ell }^{1}[g_{01},g_{02},h_{0},g_{04},\cdots,g_{0n}](g_{1\ell })  \\[4mm]
\hskip0.5cm-\delta_{b_\ell,1} \mathcal R_{\ell }^{1}[ g_{01},g_{11},g_{21},\cdots, h_0,\cdots,h_q,\cdots ,g_{0n},g_{1n},g_{2n}](g_{0\ell })\\[4mm]
\hskip0.5cm- G_{\ell }(g_{01},g_{02} ,h_{0},g_{04},\cdots, g_{0n} ,v) \bigg)dv=0,\\[4mm]
 \displaystyle\int_{V}\bigg(\partial_t h_{0} + v\cdot \nabla_{x} h_1- G_{3}(g_{01}, g_{02},h_{0},g_{04},\cdots, g_{0n} ,v) \bigg) dv=0,
\end{array}
 \right.
\end{equation}
Using (\ref{assH3}), (\ref{TZ}), (\ref{R12}) and (\ref{nS1})-(\ref{ph}), denoting by $<\cdot>$ the integral with respect to the variables $v$, shows that the system (\ref{cc1}) can be  rewritten  as follows:
\begin{equation}\label{aa}
\left\{
 \begin{array}{llll}
  \displaystyle \partial_t u_\ell+\left\langle v \cdot \nabla_{x}(\mathcal T_{\ell}^{0})^{-1}(v M_{\ell}\cdot\nabla_{x} u_\ell)\right\rangle \\[4mm]
\hskip.5cm -\delta_{b_\ell,1}\,\left\langle (\mathcal T_{\ell}^{0})^{-1}(\mathcal T_{\ell}^{1}[M_{1} u_1,\cdots, M_{n} u_n](M_{\ell} u_\ell))\right\rangle\\[4mm]
\hskip.5cm  - \left\langle G_{\ell}(M_{1} u_1,\cdots,  M_{n} u_n,v)\right\rangle = 0,\\[4mm]
\displaystyle \partial_t u_3 +\delta_{q,1}\,\left\langle v\cdot\nabla_x\mathcal T_{3}^{-1}(v M_{3}\cdot\nabla_{x} u_3)\right\rangle\\[4mm]
\hskip.5cm    - \left\langle G_{3}(M_{1} u_1,\cdots, M_{n} u_n ,v)\right\rangle  = 0,\\[4mm]
\end{array} \right.
\end{equation}

As $\mathcal{T}_{\ell}^{0}(\ell=1,2,4,\cdots,n),$ $\mathcal{T}_3$  are self-adjoint operators in $L^{2}\left(V,{dv \over M_\ell(v)}\right)$, $L^{2}\left(V ,{dv \over M_3(v)}\right)$ one has the following:
$$
\left\langle v. \nabla_{x} (\mathcal{T}_{\ell}^{0})^{-1}(v M_\ell \cdot \nabla_{x} u_\ell)\right\rangle=\div_x\bigg( \langle v \otimes
\theta_\ell(v)\rangle\cdot \nabla_{x} u_\ell\bigg),
$$
$$
\left\langle v\cdot \nabla_{x} \mathcal{T}_3^{-1}(v M_3 \cdot \nabla_{x} u_3)\right\rangle= \div_x\bigg( \langle v \otimes
\theta_3(v)\rangle\cdot \nabla_{x} u_3\bigg),
$$
and
$$
 \left\langle v\cdot \nabla_{x} (\mathcal{T}_{\ell}^{0})^{-1}(\mathcal{T}_\ell^1[M_{1} u_1,\cdots, M_{n} u_n](M_\ell u_\ell))\right\rangle  = \div_x \left\langle \frac{\theta_\ell(v)}{M_\ell(v)}\,u_\ell\mathcal{T}_\ell^1[M_{1} u_1,\cdots, M_{n} u_n](M_\ell)\right \rangle,
$$
 where $\theta_1,  \theta_2,$  $\cdots, ~\theta_n$ are given in Lemma 2.

Therefore, the macroscopic model (\ref{aa}) can be written as follows:
\begin{equation}\label{mM22}
\begin{cases}
\displaystyle{\partial_t u_1 = \vf_1\left(u_1, \cdots, u_n, \Delta u_1, \nabla u_1,\cdots,\nabla u_n \right) + \psi_1(u_1,\cdots,u_n)},\\[4mm]
\displaystyle{\partial_t u_2 = \vf_2\left(u_1, \cdots, u_n, \Delta u_2, \nabla u_1,\cdots,\nabla u_n \right) + \psi_2(u_1,\cdots,u_n)},\\[4mm]
\displaystyle{\partial_t u_3 = \vf_3\left(u_3, \Delta u_3, \nabla u_3 \right) + \psi_3(u_1,\cdots,u_n)},\\[4mm]
\displaystyle{\partial_t u_4= \vf_4\left(u_1, \cdots, u_n, \Delta u_4, \nabla u_1,\cdots,\nabla u_n \right) + \psi_4(u_1,\cdots,u_n)},\\[4mm]
\displaystyle{ \ldots ,}\\[2mm]
\displaystyle{\partial_t u_n = \vf_n\left(u_1, \cdots, u_n, \Delta u_n, \nabla u_1,\cdots,\nabla u_n\right) + \psi_n(u_1,\cdots,u_n)},\\[4mm]
\end{cases}
\end{equation}
where  $\vf_3,$  $\vf_\ell(\ell=1,2,4,\cdots,n),$  and $\psi_i(u_1,\cdots,u_n)(i=1,2,\cdots, n)$ are given by
\begin{equation} \label{w302}
\vf_\ell \left(u_1, \cdots, u_n, \Delta u_\ell, \nabla u_\ell,\cdots,\nabla u_n \right) =  \div_{x} \, (D_{\ell}\cdot \nabla_x u_\ell-\delta_{b_\ell,1}\,  u_\ell\, \alpha_\ell(u_1,\cdots,u_n)),
\end{equation}
\begin{equation} \label{w312}
\vf _3\left(u_3, \Delta u_3, \nabla u_3 \right)= \delta_{q,1} \div_{x} \, (D_{3}\cdot \nabla_x u_3),
\end{equation}
and
\begin{equation} \label{w322}
\psi_i(u_1,\cdots,u_n)= \int_V G_{i}(M_1 u_1, \cdots, M_nu_n ,v) dv.
\end{equation}
While  $D_{i}$  and $\alpha_\ell$ are given, respectively, by
\begin{equation}\label{di}
 D_{i} =- \int_V v \otimes \theta_i(v) dv,
\end{equation}
\begin{equation} \label{w1}
 \alpha_\ell(u_1,\cdots,u_n)= - \int_V  {\theta_\ell(v)\over M_\ell(v)}  \mathcal{T}_{\ell}^{1}[M_{1} u_1,\cdots,M_{n} u_n](M_{\ell})dv.
\end{equation}

\section{Micro-macro derivation  applied to selected case studies}\label{Sec:4}

This section considers the application of the general method proposed in Section 3 to the micro-macro derivation of the case studies selected in Section 2. The derivation considers, as we shall see, well defined assumptions for the interaction terms at the microscopic scale.

\subsection{SIR model in a Keller-Segel system}\label{Sec:4.1}

Let us  consider the following kernels:
\begin{equation}\label{TO}
T_\ell^0(v,v^{*})=\sigma_\ell M_\ell(v), \quad T_3(v,v^{*})=\sigma_3 M_3(v),
\end{equation}
with $\sigma_3,\,\sigma_\ell(\ell=1,2,4,\cdots, n)>0$.

Hence, the leading turning operators $\mathcal{T}^0_\ell,$ and $\mathcal{T}_3$  can be viewed as relaxation operators:
\begin{equation}\label{relaxation1}
 \mathcal{T}_\ell^{0}(g)= -\sigma_\ell \Big(g- M_\ell \langle g\rangle\Big),\hskip1cm \mathcal{T}_3(g)= -\sigma_3\Big(g-  M_3 \langle g\rangle\Big).
 \end{equation}
Moreover, $\theta_i(i=1,2,\cdots,n)$,  are given by
\[
\theta_i(v)= -\frac{1}{ \sigma_i}\,  v M_i (v),
\]
while $\alpha_\ell$, are defined by (\ref{w1}), are computed as follows:
\begin{equation}
\alpha_\ell(u_1,\cdots ,u_n)=  \frac{1}{\sigma_\ell} \int_V v \mathcal{T}_\ell^1[M_1 u_1,\cdots ,M_n u_n](M_\ell(v) )dv.\label{alpha}
\end{equation}

The diffusion tensors $D_{i}$  are given by
\begin{equation}\label{df}
D_{i}= \frac{1}{\sigma_i}\int_V v \otimes v M_i(v) dv,
 \end{equation}

while $\psi_i$  are still given by (\ref{w322}).

\subsection{Oncolytic virus model in a cross diffusion system}\label{Sec:4.2}

\vskip.2cm \noindent $\bullet$ We consider the case where the set for velocity is the sphere of radius $R > 0$, $V =R\, \mathcal{S}^{d-1}$. Let us also consider that $n=4, q=2 ~~and  ~~b_1=b_2=b_4=1$, and  the following choice:
\begin{equation}\label{Q}
T_\ell^1[f_1,f_2,f_3]=- \frac{\xi_\ell\sigma_\ell\,d}{|V|^2\,R^2}\frac{v^{*}\cdot \nabla_{x} \frac{f_3}{M_3}}{M_\ell(v^{*})},\quad \ell=1,2,4.
\end{equation}

Then $\mathcal T_\ell^1$ satisfies  (\ref{assH3}),  (\ref{R1})-(\ref{R12}), and leads to the following:
\begin{equation*}
\mathcal{T}^1_{\ell}[M_1u_1,M_2u_2,M_3u_3](M_{\ell})= \frac{\xi_{\ell}\sigma_\ell\,d}{|V|\,R^2}\,v\cdot \nabla_x u_3.
\end{equation*}
Finally,  $\alpha_\ell(u_1,u_2,u_3)$, defined in \eqref{alpha}, is given by $\alpha_\ell(u_1,u_2,u_3)= \chi_\ell\cdot \nabla_x u_3$, where the chemotactic sensitivity $\chi_\ell$ is given by the matrix
\begin{equation}\label{C}
\chi_\ell=\frac{\xi_{\ell}\,d}{|V|\,R^2}\int_V v\otimes v dv=\xi_{\ell}\,I.
\end{equation}

Therefore, the macroscopic model (\ref{mM22}) can be written as follows:

\begin{equation}\label{mM220}
\left\{
\begin{array}{l}
 \partial_t u_1=D_{1}\,\Delta u_1- \xi_{1}\,\nabla \cdot(u_1\,\nabla u_3)+\psi_1(u_1,u_2,u_3,u_4),  \\[4mm]
 \partial_t u_2=D_{2}\,\Delta u_2- \xi_{2}\,\nabla \cdot(u_2\,\nabla u_3)+\psi_2(u_1,u_2,u_3,u_4),    \\[4mm]
 \partial_t u_3= \psi_3(u_1,u_2,u_3,u_4),  \\[4mm]
 \partial_t u_4=D_{4}\,\Delta u_4- \xi_{4}\,\nabla \cdot(u_4\,\nabla u_3)+\psi_4(u_1,u_2,u_3,u_4).
\end{array} \right. \end{equation}

The role of the terms $\psi_i(u_1,u_2,u_3,u_4)(i=1,2,3,4),$  in
(\ref{w1}) consists in modeling the interaction between the four quantities of the mixture.  For example,  by choosing:
\begin{equation} \label{KK}
 G_1(f_1, f_2,f_3,f_4,v)=\frac{\mu_1}{|V|} \frac{f_1}{M_1}\big(1-\big(\frac{f_1}{M_1}\big)^{r}\big) -\frac{\rho}{|V|}\frac{f_1}{M_1} \frac{f_4}{M_4}\frac{1}{k_{u_{1}}+\theta \frac{f_1}{M_1}},
\end{equation}

\begin{equation} \label{KK1}
G_2(f_1, f_2,f_3,f_4,v)=\frac{\rho}{|V|}\frac{f_1}{M_1} \frac{f_4}{M_4}\frac{1}{k_{u_{1}}+\theta \frac{f_1}{M_1}} -\frac{\delta_2}{|V|}\frac{f_2}{M_2},
\end{equation}

\begin{equation} \label{KK2}
G_3(f_1, f_2,f_3,f_4,v)=-\frac{1}{|V|} \frac{f_3}{M_3}\big(\alpha_{1}\frac{f_1}{M_1}+\alpha_{2}\frac{f_2}{M_2}\big)
+\frac{\mu_3}{|V|} \frac{f_3}{M_3}\big(1-\frac{f_3}{M_3}\big),
\end{equation}
and
\begin{equation} \label{KK3}
G_{4}(f_1, f_2,f_3,f_4,v)=\frac{\beta}{|V|} \,\frac{f_2}{M_2}-\frac{\delta_4}{|V|} \,\frac{f_4}{M_4} -\frac{\rho}{|V|}\frac{f_1}{M_1} \frac{f_4}{M_4}\frac{1}{k_{u_{1}}+\theta \frac{f_1}{M_1}}.
\end{equation}
Hence:
\begin{equation*} \label{KK0}
      \psi_1(u_1,u_2,u_3,u_4)=\mu_1u_1(1-u_1^r) - \rho\, \frac{u_1 u_4 }{k_{u_{1}}+\theta\,u_1},
\end{equation*}
\begin{equation*} \label{KK2}
 \psi_2(u_1,u_2,u_3,u_4) =\rho \,\frac{u_1 u_4 }{k_{u_{1}}+\theta\,u_1}- \delta_2 u_2,
\end{equation*}
\begin{equation*} \label{KK00}
 \psi_3(u_1,u_2,u_3,u_4) = -u_3\left( \alpha_1 u_1 +  \alpha_2 u_2\right) + \mu_{3} u_3(1-u_3),
\end{equation*}
\begin{equation*} \label{KK20}
 \psi_4(u_1,u_2,u_3,u_4) =\beta u_2 -\delta_4 u_4- \rho \,\frac{u_1 u_4 }{k_{u_{1}}+\theta\,u_1}.
\end{equation*}
Then
\begin{equation}\label{oncolytic1}
\begin{cases}
\displaystyle{ \partial_t u_1 = D_1 \Delta u_1 - \xi_1 \nabla \cdot \left(u_1\nabla u_3\right) +\mu_1u_1(1-u_1^r) - \frac{\rho\, u_1 u_4 }{k_{u_{1}}+\theta\,u_1},}\\[4mm]
\displaystyle{ \partial_t u_2 = D_2 \Delta u_2 - \xi_2 \nabla \cdot \left(u_2 \nabla  u_3\right) + \frac{\rho \,u_1 u_4 }{k_{u_{1}}+\theta\,u_1}- \delta_2 u_2},\\[4mm]
\displaystyle{ \partial_t u_3 =-u_3\left( \alpha_1 u_1 +  \alpha_2 u_2\right) + \mu_{3} u_3(1-u_3)},\\[4mm]
\displaystyle{ \partial_t u_4 = D_4 \Delta u_4  - \xi_4 \nabla \cdot \left(u_4 \nabla u_3\right) +  \beta u_2 -\delta_4 u_4- \frac{\rho\, u_1 u_4 }{k_{u_{1}}+\theta\,u_1} }.
\end{cases}
\end{equation}

\vskip.2cm  \noindent $\bullet$  We assume that the $n=3, b_1=1, $  and $ q= b_2=2$. Then  the macroscopic model (\ref{mM22}) writes:
\begin{equation}\label{mM2021}
\left\{
\begin{array}{l}
\partial_t u_1=D_1 \Delta u_1- \div_{x} \, ( u_1\, \alpha_1(u_1,u_2,u_3))+ \psi_1(u_1,u_2,u_3),  \\[4mm]
 \partial_t u_2=D_2\Delta u_2+ \psi_2(u_1,u_2,u_3),  \\[4mm]
 \partial_t u_3= \psi_3(u_1,u_2,u_3).
\end{array} \right.
\end{equation}
In the following, we consider the kernels $T_1^0, T_2^0$ given by (\ref{TO}), $M_1(v)=M_2(v)=\frac{1}{|V|}$  and
$\sigma_1=\frac{R^2}{d},$ then
\begin{equation}\label{df}
D_{1}=I,\quad D_{2}= \frac{1}{\sigma}\,I
 \end{equation}
where $\sigma=\frac{d\,\sigma_2}{R^2}$
and $T_1^1[f_1,f_2,f_3]$ given by
\begin{equation}\label{QQ}
T_1^1[f_1,f_2,f_3]=- \frac{1}{|V|}(\chi\,v^{*}\cdot \nabla_{x}f_2+\xi\,v^{*}\cdot \nabla_{x}f_3),
\end{equation}
while the interaction terms $G_i(i=1,2,3)$  in the following form:
\begin{equation} \label{KKs}
G_1(f_1, f_2,f_3,v)=\frac{\mu}{|V|} \,\frac{f_1}{M_1}(r-\frac{f_1}{M_1}-\frac{f_3}{M_3}) ,
\end{equation}
\begin{equation} \label{KK1s}
G_2(f_1, f_2,f_3,v)=-\frac{1}{\sigma\,|V|}(\frac{f_2}{M_2}-\frac{f_1}{M_1}),
\end{equation}
\begin{equation} \label{KK2s}
G_3(f_1, f_2,f_3,v)=-\frac{1}{|V|} \,\frac{f_2}{M_2} \,\frac{f_3}{M_3} +\frac{\eta}{|V|} \,\frac{f_3}{M_3}(1-\frac{f_1}{M_1}-\frac{f_3}{M_3}) .
\end{equation}
Therefore, the macroscopic model (\ref{mM2021}) writes:
\begin{equation}\label{1}
\begin{cases}
  \displaystyle \partial_t u_1 =  \Delta u_1 - \chi \nabla\cdot (u_1\nabla u_2) - \xi \nabla\cdot (u_1\nabla u_3)+\mu\,u_1(r-u_1-u_3), \\[4mm]
  \displaystyle  \partial_t u_2 = \frac{1}{\sigma} \Delta u_2 -\frac{1}{\sigma}(u_2-u_1),   \\[4mm]
    \displaystyle \partial_t u_3 =-u_2\,u_3 +	\eta\,u_3(1-u_1-u_3).
     \end{cases}
    \end{equation}

\vskip.2cm  \noindent $\bullet$  We assume that the $n=3, b_1=b_2=q=1$. Then  the macroscopic model (\ref{mM22}) writes:
\begin{equation}\label{m1}
\left\{
\begin{array}{l}
\partial_t u_1=D_1 \Delta u_1- \div_{x} \, ( u_1\, \alpha_1(u_1,u_2,u_3))+ \psi_1(u_1,u_2,u_3),  \\[4mm]
 \partial_t u_2=D_2 \Delta u_2- \div_{x} \, ( u_2\, \alpha_2(u_1,u_2,u_3))+ \psi_2(u_1,u_2,u_3),   \\[4mm]
 \partial_t u_3= D_ 3\Delta u_3+ \psi_3(u_1,u_2,u_3).
\end{array} \right.
\end{equation}

In the following, we consider the kernels $T_1^0, T_2^0, T_3$ given by (\ref{TO}), $M_1(v)=M_2(v)=\frac{1}{|V|}$  and
$\sigma_1=\sigma_2=\frac{R^2}{d},$
and $T_1^1[f_1,f_2,f_3],  T_2^1[f_1,f_2,f_3],$ given by
$$
T_1^1[f_1,f_2,f_3]=- \frac{1}{|V|}\xi_1\,v^{*}\cdot \nabla_{x}f_3,\quad T_2^1[f_1,f_2,f_3]=- \frac{1}{|V|}\xi_2\,v^{*}\cdot \nabla_{x}f_1,
$$

while the interaction terms $G_3$  in the following form:
\begin{equation} \label{KK2s}
G_3(f_1, f_2,f_3,v)=-\frac{\lambda}{|V|} \,(\frac{f_1}{M_1} +\frac{f_2}{M_2})\frac{f_3}{M_3}-\frac{\mu}{|V|} \,\frac{f_3}{M_3}+\frac{1}{|V|}\,r(f_3) ,
\end{equation}
where $r$ satisfying the following
$$r(f+\varepsilon\,g)=r(f)+ O(\varepsilon),~~\forall f,g.$$

Therefore, the macroscopic model (\ref{m1}) writes:
\begin{equation}\label{1}
\begin{cases}
  \displaystyle \partial_t u_1 =  \Delta u_1 - \xi_1 \nabla\cdot (u_1\nabla u_3) +\psi_1(u_1,u_2,u_3), \\[4mm]
  \displaystyle  \partial_t u_2 =  \Delta u_2 - \xi_2 \nabla\cdot (u_2\nabla u_1) +\psi_2(u_1,u_2,u_3),   \\[4mm]
    \displaystyle \partial_t u_3 = D_ 3\Delta u_3-\lambda\,(u_1+u_2)\,u_3-\mu\,u_3+r(u_3).
     \end{cases}
    \end{equation}

\vskip.2cm  \noindent $\bullet$  We assume that the $n=2, b_1=1, b_2=2$. Then  the macroscopic model (\ref{mM22}) writes:
\begin{equation}\label{Q1Q}
\left\{
\begin{array}{l}
\partial_t u_1=D_1 \Delta u_1- \div_{x} \, ( u_1\, \alpha_1(u_1,u_2))+ \psi_1(u_1,u_2),  \\[4mm]
 \partial_t u_2=D_2 \Delta u_2+ \psi_2(u_1,u_2).
\end{array} \right.
\end{equation}

In the following, we consider the kernels $T_1^0, T_2^0$ given by (\ref{TO}), $$M_2(v)=\frac{1}{|V|},~~  and~~
\sigma_1=|\beta_1(u_1) \cdot\beta_2(u_2)|,$$ \\
 where $\beta_i(u_i)$(i=1,2) is a  vector valued function, and $T_1^1[f_1,f_2],$ given by
$$
T_1^1[f_1,f_2]=- \frac{\sigma_1\,d}{|V|^2\,R^2\,M_1(v^{*})}\,\frac{v^{*}\cdot \nabla_{x}\,f_2}{f_2},
$$

while the interaction terms $G_1, G_2$  in the following form:
\begin{equation} \label{KK2=s}
G_1(f_1, f_2,v)=-\frac{1}{|V|} \,\frac{f_1}{M_1} \,\frac{f_2}{M_2}+\frac{1}{|V|}\,H(f_1,f_2) .
\end{equation}
and
\begin{equation} \label{KK21=s}
G_2(f_1, f_2,v)=-\frac{1}{|V|} \,\frac{f_2}{M_2}+\frac{1}{|V|} \,\frac{f_1}{M_1} \,\frac{f_2}{M_2}+\frac{1}{|V|}\,K(f_1,f_2) .
\end{equation}
with $H, K$ satisfying the following assumptions:
$$H(f+\varepsilon\,g, \hat{f}+\varepsilon\,\hat{g})=H(f,\hat{f})+ O(\varepsilon),~K(f+\varepsilon\,g,\hat{f}+\varepsilon\,\hat{g})=K(f,\hat{f})+ O(\varepsilon),~\forall f,g,\hat{f},\hat{g}.$$

Therefore, the macroscopic model (\ref{Q1Q}) writes:
\begin{equation}\label{155}
\begin{cases}
  \displaystyle \partial_t u_1 =  \nabla\cdot \big(D(u_1,u_2)\,\nabla u_1 -  \frac{u_1}{u_2}\nabla u_2\big) -u_1\,u_2+H(u_1,u_2), \\[4mm]
  \displaystyle  \partial_t u_2 =  \eta\,\Delta u_2  -u_2 +u_1\,u_2+K(u_1,u_2),
     \end{cases}
    \end{equation}
where $$\eta=\frac{R^2}{\sigma_1\,d},~~ D(u_1,u_2)=\frac{1}{|\beta_1(u_1) \cdot\beta_2(u_2)|}\int_V v \otimes v M_1(v) dv.$$

\section{A forward look to research perspectives}\label{Sec:5}

This paper has shown how macroscopic models, which describe the dynamics of two interacting systems, can be derived at the macroscopic scale by a method somewhat inspired to the sixth Hilbert's problem. In details, we have considered the interaction between a reaction-diffusion system and models of virus dynamics described at the level of population dynamics.

The study of these systems is important as it can be specifically referred to the in-host dynamics of the SARS-CoV-2 virus, see~\cite{[BBO22],[KFPS22],[RCH22]}. A further development might be focused on the derivation of models of diffusion of epidemics  in a territory, where diffusion and reaction can be referred to networks and transportation dynamics, see~\cite{[ADKV21],[BBC20],[BDP21],[BLPZ22]}.

Technical developments  by studying models with anomalous diffusions this would be in the line of~\cite{[BC19],[BC22]} referring to models in complex environments~\cite{[BOSTW22]} with the challenging objective of understanding how a physically consistent interpretation of interactions at the microscopic scale can generate new models of nonlinear diffusion. These models have generated interesting analytic problems~\cite{[BW17]}.

A more general perspective would be  the development of the  approach in Section 3 to the micro-macro derivation of models that refer to mathematical structures analogous to that proposed in Section 2, but where the source term is modeled by different types of mathematical structures.
An immediate example is the model of reaction diffusion coupled with a Navier Stokes fluid introduced in~\cite{[Lorz12]}, followed by~\cite{[BBC16]} with the micro-macro derivation based on a first order perturbation. Further, one may consider the coupling of macroscopic equations with stochastic dynamical systems whose dynamics is described by methods of the kinetic theory of active particles. This perspective applies to self-propelled particles such as human crowds~\cite{[BGO19],[BGQR22]}. These models, as reviewed in~\cite{[BBDG21],[BCB23]}, can take into account various aspects of the heterogeneity and evolutive features of living systems.


\end{document}